\documentclass[12pt]{amsart}
\usepackage[latin1]{inputenc}
\usepackage{graphicx}
\usepackage{amssymb}
\usepackage{colordvi}
\usepackage{graphicx}
\usepackage{url}
\usepackage{enumerate}

\usepackage{vmargin}
\setpapersize{USletter}
\setmargrb{2.5cm}{2.5cm}{2.5cm}{2.5cm} % --- sets all four margins.  Cool.

\numberwithin{equation}{section}
\newtheorem{theorem}{Theorem}
\newtheorem{proposition}[theorem]{Proposition}
\newtheorem{lemma}[theorem]{Lemma}

\theoremstyle{remark}

\newtheorem{definition}[theorem]{Definition}

\newcommand{\C}{{\mathbb C}}
\newcommand{\G}{{\mathbb G}}
\newcommand{\K}{{\mathbb K}}
\renewcommand{\P}{{\mathbb P}}

\newcommand{\slC}{\mathfrak{sl}_2\C}
\newcommand{\bv}{{\bf v}}
\newcommand{\be}{{\bf e}}

\newcommand{\calK}{{\mathcal K}}

\newcommand{\calS}{{\mathcal S}}
\newcommand{\calA}{{\mathcal A}}
\newcommand{\defcolor}[1]{\RoyalBlue{#1}}
\newcommand{\demph}[1]{\defcolor{{\sl #1}}}

\def\Kosfi#1#2#3#4#5#6#7#8#9{\begin{picture}(58,25)(-2.9,-2.5)
  \put(-3.4,-3){\includegraphics[height=24pt]{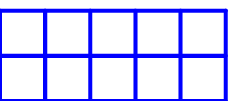}}
  \put( 0,12){\small#1}\put(12,12){\small#2}
  \put(24,12){\small#3}\put(36,12){\small#4}
  \put(48,12){\small#5}
  \put( 0, 0){\small#6}\put(12, 0){\small#7}
  \put(24, 0){\small#8}\put(36, 0){\small#9}
  \put(48, 0){\small#9}
 \end{picture}}
\def\Kosfo#1#2#3#4#5#6#7{\begin{picture}(47,25)(-2.9,-2.5)
  \put(-3.4,-3){\includegraphics[height=24pt]{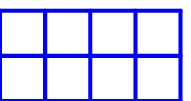}}  
  \put( 0,12){\small#1}\put(12,12){\small#2}
  \put(24,12){\small#3}\put(36,12){\small#4}
  \put( 0, 0){\small#5}\put(12, 0){\small#6}
  \put(24, 0){\small#7}\put(36, 0){\small#7}
 \end{picture}}
\title{An inequality of Kostka numbers and Galois groups of Schubert problems}
\author[Brooks]{Christopher J.\ Brooks}
\address{Christopher J.\ Brooks\\
         Department of Mathematics\\
         Texas A\&M University\\
         College Station\\
         Texas \ 77843\\
         USA}
\email{cbrooks90@neo.tamu.edu}
\author[Mart\'in del Campo]{Abraham Mart\'in del Campo}
\address{Abraham Mart\'in del Campo\\
         Department of Mathematics\\
         Texas A\&M University\\
         College Station\\
         Texas \ 77843\\
         USA}
\email{asanchez@math.tamu.edu}
\urladdr{\url{http://www.math.tamu.edu/~asanchez}}
\author[Sottile]{Frank Sottile}
\address{Frank Sottile \\
         Department of Mathematics\\
         Texas A\&M University\\
         College Station\\
         Texas \ 77843\\
         USA}
\email{sottile@math.tamu.edu}
\urladdr{\url{http://www.math.tamu.edu/~sottile}}
\thanks{Research supported in part by NSF grant DMS-915211 and the Institut
  Mittag-Leffler}
\begin{document}

\begin{abstract}
%\paragraph{Abstract.}
We show that the Galois group of any Schubert problem involving lines in projective space
 contains the alternating group. 
 Using a criterion of Vakil and a special position argument due to Schubert,
 this follows from a particular inequality among Kostka
 numbers of two-rowed tableaux.
 In most cases, an easy combinatorial injection proves the inequality.
 For the remaining cases, we use that these Kostka numbers appear in tensor product
 decompositions of $\slC$-modules.
 Interpreting the tensor product as the action of certain commuting Toeplitz matrices and using
 a spectral analysis and Fourier series rewrites the inequality as the positivity of an integral. 
 We establish the inequality by estimating this integral.
   
% \paragraph{R\'esum\'e.}
%%  texte en fran\c{c}ais
%On montre que le groupe de Galois de tout probl\`{e}me de Schubert concernant des 
%droites dans l'espace projective contient le groupe altern\'{e}.  On utilisant un crit\`{e}re de Vakil et 
%l'argument de position sp\'{e}ciale due \`{a} Schubert, ce r\'{e}sultat se d\'{e}duit d'une in\'{e}galit\'{e} particuli\`{e}re des nombres de Kostka des tableaux ayant deux rang\'{e}es. 
%Dans la plus part des cas, une injection combinatoriale facile montre l'in\'{e}galit\'{e}. Pour les cas restant, on utilise le fait que ces nombres de Kostka apparaissent dans la d\'{e}composition en produit tensoriel des $\slC$-modules.
%En interpr\'{e}tant le produit tensoriel comme l'action de certaines matrices de Toeplitz commutantent entre elles, et en utilisant de l'analyse spectrale et les s\'{e}ries de Fourier, on re\'{e}crit l'in\'{e}galit\'{e} comme la positivit\'{e}e d'une int\`{e}grale. L'in\'egalit\'e sera \'etablie en estimant cette int\`egrale.
\end{abstract}

\maketitle

%%%%%%%%%%%%%%%%%%%%%%%%%%%%%%%%%%%%%%%%%%%%%%%%%%%%%%%%%%%%%%%%%%%%%%%%%%%%%%
%
%\input{intro_11-11-28}
%%%%%%%%%%%%%%%%%%%%%%%%%%%%%%
%
\section*{Introduction}
The Schubert calculus of enumerative geometry~\cite{KL72} is a method to compute the
number of solutions to \demph{Schubert problems}, a class of geometric
problems involving linear subspaces. 
One can reduce the enumeration to combinatorics; for example, the number of solutions
to a Schubert problem involving lines is a Kostka number for a rectangular partition
with two parts.

\begin{figure}[htb]
\[
  \begin{picture}(314,192)
   \put(3,0){\includegraphics[height=6.6cm]{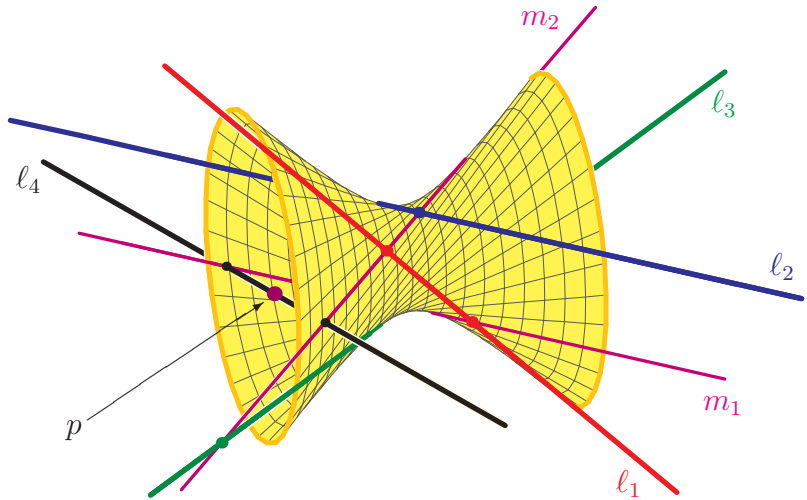}}
   \put(230,  3){$\Red{\ell_1}$}
   \put(288, 85){$\Blue{\ell_2}$}
   \put(266,147){$\Green{\ell_3}$}
   \put(  3,118){$\ell_4$}
   \put(263, 34){$\Magenta{m_1}$}
   \put(194,180){$\Magenta{m_2}$}

   \put( 22, 25){$p$}\put(30,30.5){\vector(3,2){66.6}}
  \end{picture}
\]
\caption{The two lines meeting four lines in space.\label{F:four_lines}}
\end{figure}

A prototypical Schubert problem is the classical problem of four lines, which asks for
the number of lines in space that meet four given lines. 
To answer this, note that three general lines $\Red{\ell_1},\Blue{\ell_2}$, and $\Green{\ell_3}$
lie on a unique doubly-ruled hyperboloid, shown in Figure~\ref{F:four_lines}. 
%%%%%%%%%%%%%%%%%%%%%%%%%%%%%%%%%%%%%%%%%%%%%%%%%%%%%%%%%%%%%%%%%%%%%%%%%%%%%%%%%%%%
%%%%%%%%%%%%%%%%%%%%%%%%%%%%%%%%%%%%%%%%%%%%%%%%%%%%%%%%%%%%%%%%%%%%%%%%%%%%%%%%%%%%
These three lines lie in one ruling, while the second ruling consists of
the lines meeting the given three lines.
The fourth line $\ell_4$ meets the hyperboloid in two points.
Through each of these points there is a line in the second
ruling, and these are the two lines $\Magenta{m_1}$ and $\Magenta{m_2}$ meeting our four
given lines. 
In terms of Kostka numbers, the problem of four lines reduces to counting the number of
tableaux of shape $\lambda = (2,2)$ with content $(\Red{1},\Blue{1},\Green{1},1)$. 
There are two such tableaux:
%%%%%%%%%%%%%%%%%%%%%%%%%%%%%%%%%%%%%%%%%%%%%%%%%%%%%%%%%%%%%%%%%%%%%%%%%%%%%%%%%%%%%
\[
\hspace{15pt}
\begin{picture}(56,27)(-3.4,-3)
  \put(-3.4,-3){\includegraphics{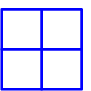}}
  \put(0,11.5){\small  \Red{1}}\put(11.5,11.5){\small \Blue{2}}
  \put(0, 0){\small \Green{3}}\put(11.5, 0){\small 4} 
\end{picture}
\hspace{-10pt}
\begin{picture}(56,27)(-3.4,-3)
  \put(-3.4,-3){\includegraphics{22.blue.eps}}
  \put(0,11.5){\small  \Red{1}}\put(11.5,11.5){\small \Green{3}}
  \put(0, 0){\small \Blue{2}}\put(11.5, 0){\small 4} 
\end{picture}
\]
%%%%%%%%%%%%%%%%%%%%%%%%%%%%%%%%%%%%%%%%%%%%%%%%%%%%%%%%%%%%%%%%%%%%%%%%%%%%%%%%%%%%%

Galois groups of enumerative problems are subtle invariants about
which very little is known.
While they were introduced by Jordan in 1870~\cite{J1870}, the modern theory began with
Harris in 1979, who showed that the algebraic Galois group is equal to a geometric monodromy
group~\cite{Ha79}. 
In general, we expect the Galois group of an enumerative problem to be the full
symmetric group and when it is not, the geometric problem possesses some intrinsic
structure. 
Harris' result gives one approach to studying the Galois group---by directly computing
monodromy.
For instance, the Galois group of the problem of four lines is the group of permutations
which are obtained by following the solutions over loops in the space of lines
$\Red{\ell_1},\Blue{\ell_2},\Green{\ell_3}, \ell_4$.
Rotating $\ell_4$ 180 degrees about the point $p$ (shown in Figure \ref{F:four_lines}) gives
a loop which interchanges the two solution lines $\Magenta{m_1}$ and $\Magenta{m_2}$, showing
that the Galois group is the full symmetric group on two letters.

Leykin and Sottile~\cite{LS09} used numerical homotopy continuation~\cite{SW05} to
compute monodromy for many {\it simple} Schubert problems, showing that in each case the
Galois group was the full symmetric group. 
(The problem of four lines is simple.)
Billey and Vakil~\cite{BV} gave an algebraic approach based on elimination theory to
compute lower bounds for Galois groups.
Vakil~\cite{Va06b} gave a combinatorial criterion, based on group theory, which can be used
to show that a Galois group contains the alternating group.
He used this and his geometric Littlewood-Richardson rule~\cite{Va06a} to show that the
Galois group was at least alternating for every Schubert problem involving lines in
projective space $\P^n$ for $n\leq 16$. 
Brooks implemented Vakil's criterion and the geometric Littlewood-Richardson
rule in {\tt python} and used it to show
that for $n\leq 40$, every Schubert problem involving lines in projective space $\P^n$ has
at least alternating Galois group. 
Our main result is the following.

%%%%%%%%%%%%%%%%%%%%%%%%%%%%%%%%%%%%%%%%%%%%%%%%%%%%%%%%%%%%%%%%%%%%%%%%
\begin{theorem}\label{Th:Only}
 The Galois group of any Schubert problem involving lines in $\P^n$ contains the 
  alternating group.
\end{theorem}
%%%%%%%%%%%%%%%%%%%%%%%%%%%%%%%%%%%%%%%%%%%%%%%%%%%%%%%%%%%%%%%%%%%%%%%%

We prove this theorem by applying Vakil's criterion to a special position argument of
Schubert, which reduces Theorem~\ref{Th:Only} to proving a certain inequality among Kostka
numbers of two-rowed tableaux.  For most problems, the inequality follows from a
combinatorial injection of Young tableaux. For the remaining problems, we work in the
representation ring of $\slC$, where these Kostka numbers also occur.
We interpret the tensor product of irreducible $\slC$-modules in terms of commuting
Toeplitz matrices.
Using the eigenvector decomposition of the Toeplitz matrices, we express these Kostka numbers as
certain trigonometric integrals. 
In this way, the inequalities of Kostka numbers become inequalities of integrals, which we
establish by estimation.

Note that the generalization of Theorem~\ref{Th:Only} to arbitrary Grassmannians is false.
Derksen found Schubert problems in the Grassmannian of $3$-planes in $\P^7$ whose Galois
groups are significantly smaller than the full symmetric group, and Vakil generalized this to 
problems in the Grassmannians of $2k{-}1$ planes in $\P^{2n-1}$ whose Galois groups are not the 
full symmetric group for every $k\geq 2$ and $n\geq 2k$~\cite[\S 3.13]{Va06b}.

%%%%%%%%%%%%%%%%%%%%%%%%
%\input{section1_11-11-28}
\section{Preliminaries}\label{S:one}

%%%%%%%%%%%%%%%%%%%%%%%%%%%%%%%%%%%%%%%%%%%%%%%%%%%%%%%%%%%%%%%%%%%%%%%%%%%%%%
%
\subsection{Schubert problems of lines}
Let $\G(1,n)$ be the Grassmannian of lines in $n$-dimensional
projective space $\P^n$, which is an algebraic manifold of dimension $2n{-}2$.
A (special) \demph{Schubert subvariety} is the set of lines $X_L$ that meet a linear subspace
$L\subset\P^n$; that is,
 \begin{equation}\label{Eq:SchubertVariety}
   \defcolor{X_L}\ :=\ \{\ell\in\G(1,n)\mid \ell\cap L\neq\varnothing\}\,.
 \end{equation}
If $\dim L= n{-}1{-}a$, then $X_L$ has codimension $a$ in $\G(1,n)$. 
A \demph{Schubert problem} asks for the lines that meet fixed
linear subspaces $L_1,\ldots,L_m$ in general position,
where  $\dim L_i = n{-}1{-}a_i$ for
$i=1,\ldots, m$ and $a_1+\cdots +a_m= 2n{-}2$.
These are the points in the intersection 
\begin{equation}\label{Eq:SchubertProblem}
   X_{L_1}\cap X_{L_2}\cap\cdots\cap X_{L_m}.
 \end{equation}
As the $L_i$ are in general position, the intersection (\ref{Eq:SchubertProblem}) is
transverse and therefore zero-dimensional.  
(Over fields of characteristic zero, transversality follows from Kleiman's Transversality
Theorem~\cite{Kl74} while in positive characteristic, it is Theorem~E in~\cite{So97}.)
We define the \demph{Schubert intersection} number $\defcolor{K(a_1,\ldots,a_m)}$ to be
the number of points in the intersection (\ref{Eq:SchubertProblem}), which does not depend
upon the choice of general $L_1,\ldots, L_m$. 
We call \defcolor{$a_\bullet:=(a_1,\ldots,a_m)$} the \demph{type} of the Schubert
problem~(\ref{Eq:SchubertProblem}). 

Note that given positive positive integers $a_\bullet=(a_1,\ldots,a_m)$ whose sum is even, 
$K(a_\bullet)$ is a Schubert intersection number in $\G(1,n(a_\bullet))$, where
$\defcolor{n(a_\bullet)}:=\frac{1}{2}(a_1+\cdots+a_m+2)$.
Henceforth, a Schubert problem will be a list $a_\bullet$ of positive integers with
even sum. 
It is \demph{valid} if $a_i\leq n(a_\bullet){-}1$ (this is forced by $\dim L_i\geq 0$).

The intersection number $K(a_\bullet)$ is a Kostka number, which is the number of Young
tableaux of shape $(n(a_\bullet){-}1,n(a_\bullet){-}1)$ and content
$(a_1,\ldots,a_m)$~\cite[p.25]{Fu97}.   
Let \defcolor{$\calK(a_\bullet)$} be the set of such tableaux.
These are two-rowed arrays of integers, each row of length
$n(a_\bullet){-}1$, such that the integers
increase weakly across each row and strictly down each column, and
there are $a_i$ occurrences of $i$ for each $i=1,\ldots,m$.
For example, here are the five Young tableaux in $\calK(2,2,1,2,3)$,
demonstrating that $K(2,2,1,2,3)=5$.
\begin{equation}\label{Eq11:tableaux}
 \raisebox{-9pt}{%
   \Kosfi{1}{1}{2}{2}{3}{4}{4}{5}{5}\quad
   \Kosfi{1}{1}{2}{2}{4}{3}{4}{5}{5}\quad
   \Kosfi{1}{1}{2}{3}{4}{2}{3}{5}{5}\quad
   \Kosfi{1}{1}{2}{4}{4}{2}{3}{5}{5}\quad
   \Kosfi{1}{1}{3}{4}{4}{2}{2}{5}{5}}
\end{equation}

%%%%%%%%%%%%%%%%%%%%%%%%%%%%%%%%%%%%%%%%%%%%%%%%%%%%%%%%%%%%%%%%%%%%%%%%%%%%%%
%
\subsection{Vakil's Criterion for Galois groups of Schubert problems}\label{S:vakil}
In \S 3.4 of~\cite{Va06b}, Vakil explains how to associate a Galois group
to a dominant map $W\to X$ of equidimensional irreducible varieties and  
establishes his criterion for the Galois group to contain the alternating group.
We discuss this for a Schubert problem $a_\bullet=(a_1,\ldots,a_m)$.
Define
\[
  \defcolor{X}\ :=\ \{(L_1,\ldots,L_m)\mid L_i\subset \P^n
  \mbox{\ is a linear space of dimension }n{-}1{-}a_i\}\,,
\]
where $n:=n(a_\bullet)$.
Consider the incidence variety, 
\[
  \defcolor{W}\ := \{(\ell,L_1,\ldots,L_m)\mid (L_1,\ldots,L_m)\in X
    \mbox{\ and\ } \ell\cap L_i\neq\varnothing\,,\ i=1,\ldots,m\}\,.
\]
The projection map $W\to\G(1,n)$ realizes $W$ as a fiber bundle over $\G(1,n)$ with
irreducible fibers.
As $\G(1,n)$ is irreducible, $W$ is irreducible.

Let $\pi\colon W\to X$ be the other projection;
its fiber over a point $(L_1,\ldots,L_m)\in X$ is
 \begin{equation}\label{Eq:fiber}
   \pi^{-1}(L_1,L_2,\ldots,L_m)\ =\ X_{L_1}\cap X_{L_2}\cap\cdots\cap X_{L_m}\,.
 \end{equation}
Thus the map $\pi\colon W\to X$ contains all Schubert problems of type
$a_\bullet$.
As the general Schubert problem is a transverse intersection containing $K(a_\bullet)$
points, $\pi$ is a dominant map of degree $K(a_\bullet)$.
Under $\pi$, the field \defcolor{$\K(X)$} of rational functions on $X$ pulls back to a
subfield of $\K(W)$, the field of rational functions on $W$, and the extension
$\K(W)/\K(X)$ has degree $K(a_\bullet)$.

\begin{definition}
  The Galois group of the Schubert problem of type $a_\bullet$, 
  \defcolor{$G(a_\bullet)$}, is the Galois group of the Galois closure of the field extension $\K(W)/\K(X)$.
\end{definition}

This Galois group $G(a_\bullet)$ is a subgroup of the symmetric group $\calS_{K(a_\bullet)}$ on $K(a_\bullet)$ letters.
We say that $G(a_\bullet)$ is \demph{at least alternating} if it contains the alternating group $\calA_{K(a_\bullet)}$.
Vakil's Criterion  %for showing that a Galois group is at least alternating 
is adapted to classical special position arguments in enumerative geometry.
First, if $Z\subset X$ is a subvariety such that $Y=\pi^{-1}(Z)\subset W$ is irreducible
and the map $Y\to Z$ has degree $K(a_\bullet)$, then $Y\to Z$ has a Galois group which is
a subgroup of $G(a_\bullet)$.
This enables us to restrict the original Schubert problem to one derived from it through
certain standard reductions.

More interesting is when $Z\subset X$ is a subvariety such that $Y=\pi^{-1}(Z)$ 
decomposes into two smaller problems, $Y=Y_1\cup Y_2$, where $Y_i\to Z$ is a Schubert
problem of type $a^{(i)}_\bullet$ for $i=1,2$.
In this situation, monodromy of $Y\to Z$ gives a subgroup $H$ of the product 
$G(a^{(1)}_\bullet)\times G(a^{(2)}_\bullet)$ which projects onto each factor and includes into $G(a_\bullet)$.
Then purely group-theoretic arguments imply the following.\medskip

%%%%%%%%%%%%%%%%%%%%%%%%%%%%%%%%%%%%%%%%%%%%%%%%%%%%%%%%%%%%%%%%%%%%%%%%%%%%%%
\noindent{\bf Vakil's Criterion.}
{\it
  If $G(a^{(1)}_\bullet)$ and $G(a^{(2)}_\bullet)$ are at least alternating, and either 
  $K(a^{(1)}_\bullet)\neq K(a^{(2)}_\bullet)$ or
  $K(a^{(1)}_\bullet)=K(a^{(2)}_\bullet)=1$; then $G(a_\bullet)$ is at least
  alternating. 
}\medskip

%%%%%%%%%%%%%%%%%%%%%%%%%%%%%%%%%%%%%%%%%%%%%%%%%%%%%%%%%%%%%%%%%%%%%%%%%%%%%%
%
%\input{section2_11-11-28}
\section{Inequalities}\label{S:two}

A Schubert problem  $a_\bullet=(a_1,\ldots,a_m)$ is \demph{reduced} if it is valid and 
if $a_i{+}a_j\leq n(a_\bullet){-}1$ for any $i<j$.
Any Schubert problem is equivalent to a reduced one:
If $a_\bullet$ is valid, but $a_{m-1}+a_m>n(a_\bullet){-}1$, then 
\[
   K(a_1,\ldots,a_m)\ =\ K(a_1,\ldots,a_{m-2},\, a_{m-1}{-}1, a_m{-}1)\,,
\]
as the intersection~(\ref{Eq:SchubertProblem}) for $a_{\bullet}$ is equal to an intersection
for $(a_1,\ldots,a_{m-2},a_{m-1}{-}1, a_m{-}1)$.
Iterating this procedure gives an equivalent reduced Schubert problem.

Schubert~\cite{Sch1886a} observed that if the linear spaces are in a special position,
then the Schubert problem decomposes into two smaller problems, which gives a (familiar)
recursion for these Kostka numbers. 
Given a reduced Schubert problem $a_\bullet=(a_1,\ldots,a_m)$, set $n:=n(a_\bullet)$.
Let $L_1,\ldots,L_m$ be linear subspaces which are in general position in $\P^n$, except that
$L_{m-1}$ and $L_m$ span a hyperplane $\Lambda:= \overline{L_{m-1},L_m}$. 
If a line $\ell$ meets both $L_{m-1}$ and $L_m$, then either it meets $L_{m-1}\cap L_{m}$
or it lies in their linear span (while also meeting both $L_{m-1}$ and $L_{m}$). 
This implies Schubert's recursion for Kostka numbers
 \begin{equation}\label{E:combinatorial_recursion}
  K(a_1,\ldots,a_m)\ =\ 
    K(a_1,\ldots,a_{m-2},\,a_{m-1}{+}a_m)\ +\ 
    K(a_1,\ldots,a_{m-2},\, a_{m-1}{-}1,a_m{-}1)\,.
 \end{equation}
Observe that if $a_\bullet$ is reduced, then both smaller problems
in~(\ref{E:combinatorial_recursion}) are valid.
An induction shows that if $a_\bullet$ is valid, then $K(a_\bullet)>0$.

For example, consider $K(2,2,1,2,3)$.
The first tableau in~(\ref{Eq11:tableaux}) has both 4s in its second row (along with its
5s), while the remaining four tableaux have last column consisting of a 4 on top of a 5.
If we replace the 5s by 4s in the first tableau and erase the last column in the remaining
four tableaux, we obtain
\[
   \Kosfi{1}{1}{2}{2}{3}{4}{4}{4}{4}\qquad\quad
   \Kosfo{1}{1}{2}{2}{3}{4}{5}\quad
   \Kosfo{1}{1}{2}{3}{2}{3}{5}\quad
   \Kosfo{1}{1}{2}{4}{2}{3}{5}\quad
   \Kosfo{1}{1}{3}{4}{2}{2}{5}\ \,
\]
which shows that $K(2,2,1,2,3)=K(2,2,1,5)+K(2,2,1,1,2)$.
We state our key lemma.
A \demph{rearrangement} of a Schubert problem $a_1,\ldots,a_m$ is simply a listing of the
integers $a_1,\ldots,a_m$ in some order.

%%%%%%%%%%%%%%%%%%%%%%%%%%%%%%%%%%%%%%%%%%%%%%%%%%%%%%%%%%%%%%%%%%%%%%%%%%%%%%
\begin{lemma}\label{L:induction}
 Every reduced Schubert problem has a rearrangement $(a_1,\ldots,a_m)$ such that either
 \begin{equation}\label{Eq:ineq}
     K(a_1,\ldots,a_{m-2},\,a_{m-1}{+}a_m)\ \neq\ 
    K(a_1,\ldots,a_{m-2},\, a_{m-1}{-}1,a_m{-}1)\,,
 \end{equation}
 and both are nonzero, or else both are equal to $1$.
\end{lemma}
%%%%%%%%%%%%%%%%%%%%%%%%%%%%%%%%%%%%%%%%%%%%%%%%%%%%%%%%%%%%%%%%%%%%%%%%%%%%%%

We use Lemma~\ref{L:induction} below to prove Theorem~\ref{Th:Only}, then we devote the rest
of the extended abstract to the proof of this Lemma.

%%%%%%%%%%%%%%%%%%%%%%%%%%%%%%%%%%%%%%%%%%%%%%%%%%%%%%%%%%%%%%%%%%%%%%%%%%%%%%
\begin{proof}[of Theorem~\ref{Th:Only}]
 We use the notation of Subsection~\ref{S:vakil} and argue by induction on $m$ and $n(a_\bullet)$.
 Assume that $a_\bullet$ is reduced and let $Z$ be the set of those
 $(L_1,\ldots,L_m)\in X$ such that $L_{m-1}$ and $L_m$ span a hyperplane.
 Then the geometric arguments given before~(\ref{E:combinatorial_recursion}) imply that
 the pullback $\pi^{-1}(Z)\to Z$ decomposes as the union of two Schubert problems, one for 
 $(a_1,\ldots,a_{m-2},\,a_{m-1}{+}a_m)$ and the other for $(a_1,\ldots,a_{m-2},\, a_{m-1}{-}1,a_m{-}1)$.
 Therefore, Lemma~\ref{L:induction} and our induction hypothesis, together with Vakil's
 criterion, imply that $G(a_\bullet)$ is at least alternating.
\end{proof}
%%%%%%%%%%%%%%%%%%%%%%%%%%%%%%%%%%%%%%%%%%%%%%%%%%%%%%%%%%%%%%%%%%%%%%%%%%%%%%

While an induction shows that the only reduced Schubert problem where the two terms
in~(\ref{Eq:ineq}) are both 1 is $(1,1,1,1)$, the inequality of Lemma~\ref{L:induction} 
is not easy to prove. This is in part because there are no 
closed formulas for the numbers $K(a_\bullet)$, except for the case $a_1=\cdots =a_{m-1} =
1$ (in which case $K(a_\bullet)$ is given by the hook-length formula).

%%%%%%%%%%%%%%%%%%%%%%%%%%%%%%%%%%%%%%%%%%%%%%%%%%%%%%%%%%%%%%%%%%%%%%%%%%%%%%
%
\subsection{Inequality of Lemma~\ref{L:induction} in most cases}

We give an injection of sets of Young tableaux to establish
Lemma~\ref{L:induction} when $a_i\neq a_j$ for some $i,j$.

%%%%%%%%%%%%%%%%%%%%%%%%%%%%%%%%%%%%%%%%%%%%%%%%%%%%%%%%%%%%%%%%%%%%%%%%%%%%%%
\begin{lemma}\label{L:unequal_inequality}
 Suppose that $(b_1,\ldots,b_m,\alpha,\beta,\gamma)$ is a reduced Schubert problem where 
 $\alpha\leq\beta\leq\gamma$ with $\alpha<\gamma$.
 Then
 \begin{equation}\label{Eq:unequal_inequality}
   K(b_1,\ldots,b_m,\,\alpha, \beta+\gamma)\ <\ 
   K(b_1,\ldots,b_m,\,\gamma, \beta+\alpha)\,.
 \end{equation}
\end{lemma}
%%%%%%%%%%%%%%%%%%%%%%%%%%%%%%%%%%%%%%%%%%%%%%%%%%%%%%%%%%%%%%%%%%%%%%%%%%%%%%

To see that this implies Lemma~\ref{L:induction} in the case when $a_i\neq a_j$, for
some $i,j$, we apply Schubert's recursion
to obtain two different expressions for $K(b_1,\ldots,b_m,\alpha,\beta,\gamma)$,
%
% \begin{eqnarray*}
 \begin{multline*}
   \qquad   K(b_1,\ldots,b_m,\,\alpha, \beta{+}\gamma)\ + \ 
   K(b_1,\ldots,b_m,\,\alpha, \beta{-}1,\gamma{-}1)\\
  \ =\    K(b_1,\ldots,b_m,\,\gamma, \beta+\alpha)\ +\ 
   K(b_1,\ldots,b_m,\,\gamma, \beta{-}1,\alpha{-}1)\,.\qquad
 \end{multline*}
% \end{eqnarray*}
%
By the inequality~(\ref{Eq:unequal_inequality}), at least one of these expressions involves
unequal terms.
Since all four terms are from valid Schubert problems, none is zero, and this implies
Lemma~\ref{L:induction} when not all $a_i$ are identical.\hfill\qed\medskip

%%%%%%%%%%%%%%%%%%%%%%%%%%%%%%%%%%%%%%%%%%%%%%%%%%%%%%%%%%%%%%%%%%%%%%%%%%%%%%
\begin{proof}[of Lemma~$\ref{L:unequal_inequality}$]
 We establish the inequality~(\ref{Eq:unequal_inequality}) via a combinatorial injection
\[
   \iota\ \colon\  \calK(b_1,\ldots,b_m,\,\alpha, \beta+\gamma)\ 
   \lhook\joinrel\relbar\joinrel\rightarrow\ 
   \calK(b_1,\ldots,b_m,\,\gamma, \beta+\alpha)\,,
\]
 which is not surjective.

 Let $T$ be a tableau in $\calK(b_1,\ldots,b_m,\,\alpha, \beta+\gamma)$ and let $A$ be its
 sub-tableau consisting of the entries $1,\ldots,m$.
 Then the skew tableau $T\setminus A$ has a bloc of $(m{+}1)$'s of length $a$ at the end
 of its first row, and its second row consists of a bloc of $(m{+}1)$'s of length $\alpha{-}a$,
 followed by a bloc of $(m{+}2)$'s of length $\beta{+}\gamma$.
 Form the tableau $\iota(T)$ by changing the last row of $T\setminus A$ to 
 a bloc of $(m+1)$'s of length $\gamma{-}a$ followed by a bloc of $(m+2)$'s of length $\beta{+}\alpha$.
 Since $a\leq\alpha<\gamma$, this map is well-defined.
\[
  T\ =\ 
  \raisebox{-13pt}{\begin{picture}(131,33)
   \put(0,0){\includegraphics{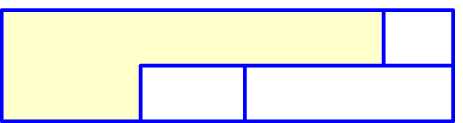}}
   \put(117,22){$a$}
   \put(45,6){$\alpha{-}a$}\put(90,6){$\beta{+}\gamma$}
   \put(17,15){$A$}
  \end{picture}}
   \ \longmapsto\ 
  \raisebox{-13pt}{\begin{picture}(131,33)
   \put(0,0){\includegraphics{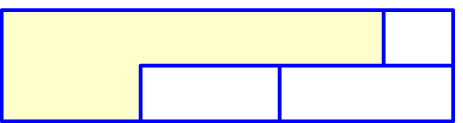}}
   \put(117,22){$a$}
   \put(50,6){$\gamma{-}a$}\put(95,6){$\beta{+}\alpha$}
   \put(17,15){$A$}
  \end{picture}}
  \ =\ \iota(T)\,.
\]
 
 To see that $\iota$ is not surjective, set
 $\defcolor{b_\bullet}:=(b_1,\ldots,b_m,\gamma{-}\alpha{-}1,\beta{-}1)$, which is a valid
 Schubert problem.
 Hence $K(b_\bullet)\neq 0$ and $\calK(b_\bullet)\neq\varnothing$.
 For any $T\in\calK(b_\bullet)$, we may add $\alpha{+}1$ columns to its end consisting of a
 $m{+}1$ above a $m{+}2$ to obtain a tableau
 $T'\in\calK(b_1,\ldots,b_m,\,\gamma,\beta+\alpha)$.
 As $T'$ has more than $\alpha$ $(m{+}1)$s in its first row, it cannot be in the image of the
 injection $\iota$, which completes the proof of the lemma.
\end{proof}
%%%%%%%%%%%%%%%%%%%%%%%%%%%%%%%%%%%%%%%%%%%%%%%%%%%%%%%%%%%%%%%%%%%%%%%%%%%%%%

%\input{section3_11-11-28}
%%%%%%%%%%%%%%%%%%%%%%%%%%%%%%%
%
%\section{Ring of Symmetric functions in two variables}\label{S:three}
\section{Kostka numbers as integrals}\label{S:three}

Kostka numbers of two-rowed tableaux appear as the coefficients in the decomposition of the
tensor products of irreducible $\slC$-modules.
Let \defcolor{$V_a$} be the irreducible module of $\slC$ with highest weight $a$. 
Given a Schubert problem $a_\bullet = (a_1,\ldots,a_m)$, the Kostka number $K(a_\bullet)$ is
the multiplicity of the trivial $\slC$-module $V_0$ in the tensor product
$V_{a_1}\otimes\cdots\otimes V_{a_m}$.   

The representation ring $R$ of $\slC$ is the free abelian group on the isomorphism classes
$[V_a]$ of irreducible modules, modulo the relations $[V_a]+ [V_b]-[V_a\oplus V_b]$. 
Setting $[V_a]\cdot [V_b]:=[V_a\otimes V_b]$ makes $R$ into a ring.
Writing $\defcolor{\be_a}:=[V_a]$, multiplication by $\be_a$ is a linear operator 
\defcolor{$M_a$} on $R$, 
\begin{equation}\label{E:operator}
M_a (\be_b)\ :=\   \be_a\cdot \be_b
 \ =\ \be_{b+a} {+} \be_{b+a-2} + \cdots + \be_{|b-a|}\,,
\end{equation}
by the Clebsch-Gordan formula.
In the basis $\{\be_a\}$, the operator $M_a$ is represented by an infinite Toeplitz matrix
with entries $0$ and $1$ given by the formula~(\ref{E:operator}).
For instance, we have
\[
M_2\ =\ \left( 
\begin{array}{cccccccc}
0 & 0 & 1 & 0 & 0 & 0 & 0 &  \\
0 & 1 & 0 & 1 & 0 & 0 & 0 &  \\
1 & 0 & 1 & 0 & 1 & 0 & 0 & \cdots \\
0 & 1 & 0 & 1 & 0 & 1 & 0 &  \\
0 & 0 & 1 & 0 & 1 & 0 & 1 &  \\ 
&  &  & \vdots &  &  &  & \ddots
\end{array}
\right),
\hspace{25pt} 
M_3\ =\ \left(
\begin{array}{ccccccccc}
0 & 0 & 0 & 1 & 0 & 0 & 0 & 0 &  \\
0 & 0 & 1 & 0 & 1 & 0 & 0 & 0 &  \\
0 & 1 & 0 & 1 & 0 & 1 & 0 & 0 & \cdots \\
1 & 0 & 1 & 0 & 1 & 0 & 1 & 0 &  \\
0 & 1 & 0 & 1 & 0 & 1 & 0 & 1 &  \\ 
&  &  & \vdots &  &  &  &  & \ddots 
\end{array}
\right).
\]
Since $R$ is a commutative ring, the operators $\{M_a\mid a\geq 0\}$ commute.
They have an easily described system of joint eigenvectors and eigenvalues, which may be
verified using the identity $2\sin\alpha\cdot\sin\beta=\cos(\alpha{-}\beta)-\cos(\alpha{+}\beta)$,  
and noting that the resulting sums are telescoping.
%%%%%%%%%
%
\begin{proposition}\label{P:eigenvectors}
 For each $0\leq \theta \leq \pi$ and integer $a\geq 0$, set
 \begin{eqnarray*}
   \demph{\bv(\theta)} & :=& 
    (\sin\theta, \sin 2\theta, \ldots, \sin(j{+}1)\theta,\ldots)^\top \ = \  \sum_j \sin (j{+}1)\theta \cdot \be_j, \\    
    \demph{\lambda_a(\theta)} &:=&  \frac{\sin (a{+}1)\theta}{\sin \theta}\,.
 \end{eqnarray*}
 Then $\bv(\theta)$ is an eigenvector of $M_a$ with eigenvalue $\lambda_a(\theta)$.
\end{proposition}

%While we do not need it, we may regard $R$ as a dense subspace of the Banach space
%$\ell_\infty$, and the $M_a$ become bounded operators on $\ell_\infty$.
%%%%%%%%%%%%%%%%%%%%%%%%%%%%%%%%%%%%%%%%%%%%%%%%%%%%%%%%%%%%%%%%%%%%%%%%%%%%%%

These eigenvectors form a complete system of eigenvectors.
%%%%%%%%%%%%%%%%%%
\begin{proposition}\label{P:eigenvector_decomposition}
For any $a=0,1,2,\ldots$, we have 
 \[
    \be_j\ =\ \frac{2}{\pi}\int_0^{\pi} \sin{(j{+}1)\theta}\ \bv(\theta)\, d\theta\,.
 \]
It follows that for any $a\geq 1$, we have
 \[
   M_a(\be_0)\ =\ \frac{2}{\pi}\int_{0}^{\pi}\lambda_a(\theta) \sin \theta \, \bv(\theta)\, 
   d\theta\,.
 \]
\end{proposition}

A consequence of Proposition~\ref{P:eigenvector_decomposition} is an integral formula for the
Kostka numbers.

%%%%%%%%%%%%%%%%%%%%%%%%%%%%%%%%%%%%%%%%%%%
\begin{theorem}
  Let $a_\bullet=(a_1,\dotsc,a_m)$ be any valid Schubert problem. Then
 \begin{equation}\label{Eq:Kostka_integral}
   K(a_\bullet)\ =\ \frac{2}{\pi}\int_{0}^{\pi} 
      \left(\prod_{i=1}^m\lambda_{a_i}(\theta)\right) \;\sin^2 \theta\, d\theta\,.
 \end{equation}
\end{theorem}
%%%%%%%%%%%%%%%%%%%%%%%%%%%%%%%%%%%%%%%%%%%

%%%%%%%%%%%%%%%%%%%%%%%%%%%%%%%%%%%%%%%%%%%%%%%%%%%
%
\subsection{Inequality of Lemma~\ref{L:induction} in the remaining case}
%%%%%%%%%%%%%%%%%%%%%%%%%%%%%%%%%%%%%%%%%%%%%%%%%%

We complete the proof of Theorem~\ref{Th:Only} by establishing the inequality in
Lemma~\ref{L:induction} for those Schubert problems not covered in Lemma~\ref{L:unequal_inequality}. For these, every condition is the same, so $a_\bullet = (a,a,\ldots,a)=:\defcolor{a^m}$.
%in the remaining cases of Schubert problems $a_\bullet$ of the
%form $(a,a,\ldots,a)=:\defcolor{a^m}$, which are not covered in Lemma~\ref{L:unequal_inequality}. 

If $a=1$, then we may use the hook-length formula.
The Kostka number $K(1^{n},b)$, where $n+b=2c$ is even, is the number of Young tableaux
of shape $(c,c-b)$, which is
\[
   K(1^{n},b)\ :=\ \frac{n!(b{+}1)}{(c{-}b)!(c{+}1)!}
\]
When $m=2n$ is even, the inequality of Lemma~\ref{L:induction} is
that $K(1^{2n-2})\neq K(1^{2n-2},2)$.
We compute
\[
  K(1^{2n-2})\ =\ \frac{(2n-2)!(1)}{n!(n+1)!}
  \qquad\mbox{and}\qquad
  K(1^{2n-2},2)\ =\ \frac{(2n-2)!(3)}{(n-2)!(n+1)!}
\]
and so 
\[
   K(1^{2n-2},2)/K(1^{2n-2})\ =\ 3\frac{n!(n{+}1)!}{(n{-}2)!(n{+}1)!}
    \ =\ 3\frac{n{-}1}{n{+}1}\ \neq\ 1\,,
\]
when $n>2$, but when $n=2$ both Kostka numbers are $1$, which proves the  inequality of
Lemma~\ref{L:induction} when each $a_i=1$. 

We now suppose that $a_\bullet=(a^{m+2})$ where $a>1$ and $m\cdot a$ is even.
Table~\ref{T:inequality_a=2} shows that when $a=2$ and $m\leq 16$, the inequality
of Lemma~\ref{L:induction} holds.
%%%%%%%%%%%%%%%%%%
\begin{table}[htdp]
\caption{The inequality \eqref{Eq:ineq} for the case $a_\bullet=(2^{m+2})$}
\label{T:inequality_a=2}
\begin{center}\begin{tabular}{|c||c|c|c|}
\hline
$m$ & $K(2^m,4)$ & $K(2^m,1,1)$ & Difference \\\hline\hline
0 & 0 & 1 & $-1$ \\\hline
1 & 0 & 1 & $-1$ \\\hline
2 & 1 & 2 & $-1$ \\\hline
3 & 2 & 4 & $-2$ \\\hline
4 & 6 & 9 & $-3$ \\\hline
5 & 15 & 21 & $-6$ \\\hline
6 & 40 & 51 & $-11$ \\\hline
7 & 105 & 127 & $-22$ \\\hline
8 & 280 & 323 & $-43$ \\\hline
9 & 750 & 835 & $-85$ \\\hline
10 & 2025 & 2188 & $-163$ \\\hline
11 & 5500 & 5798 & $-298$ \\\hline
12 & 15026 & 15511 & $-485$ \\\hline
13 & 41262 & 41835 & $-573$ \\\hline
14 & 113841 & 113634 & 207\\\hline
15 & 315420 & 310572 & 4848\\\hline
16 & 877320 & 853467 & 23853\\\hline
%17 & 2448816 & 2356779 & 92037\\\hline
%18 & 6857307 & 6536382 & 320925\\\hline
%19 & 19259046 & 18199284 & 1059762\\\hline
%20 & 54237210 & 50852019 & 3385191\\\hline
\end{tabular} 
\end{center}
\end{table}
%%%%%%%%%%%%%%%%%%
However, the sign of $K(2^m,4)-K(2^m,1,1)$ changes at $m=14$.
In fact, we have the following lemma.

%%%%%%%%%%%%%%%%%%%%%%%%%%%%%%%%%%%%%%%%%%%%%%%%%%
\begin{lemma}\label{L:inequality_a=2}
  For all $m\geq 1$, we have $K(2^m,4)\neq K(2^m,1,1)$.
  If $m<14$ then $K(2^m,4)<K(2^m,1,1)$ and if $m\geq 14$, then $K(2^m,4) >K(2^m,1,1)$.
\end{lemma}
%%%%%%%%%%%%%%%%%%%%%%%%%%%%%%%%%%%%%%%%%%%%%%%%%%

The remaining cases $a\geq 3$ have a more uniform behavior.

%%%%%%%%%%%%%%%%%%%%%%%%%%%%%%%%%%%%%%%%%%%%%%%%%%
\begin{lemma}\label{L:equal_inequality}
 For $a\geq3$ and for all $m\geq 2$ we have
 \begin{equation}\label{E:equal_inequality}
   K(a^m,\,2a)\ <\ K(a^m,(a{-}1)^2)\,.
 \end{equation}
\end{lemma}
%%%%%%%%%%%%%%%%%%%%%%%%%%%%%%%%%%%%%%%%%%%%%%%%%%

We omit the proof of Lemma~\ref{L:equal_inequality} from this extended abstract, but
include a proof of Lemma~\ref{L:inequality_a=2}.

%%%%%%%%%%%%%%%%%%%%%%%%%%%%%%%%%%%%%%%%%%%%%%%%%%%%%%%%%%%%%%%%%%%%%%%%%%%%%%
%
\subsection{Proof of Lemma~$\ref{L:inequality_a=2}$}\label{SS:a=2}

\begin{figure}[htb]

   \begin{picture}(127,153)(-7,0)
    \put(0,0){\includegraphics[height=150pt]{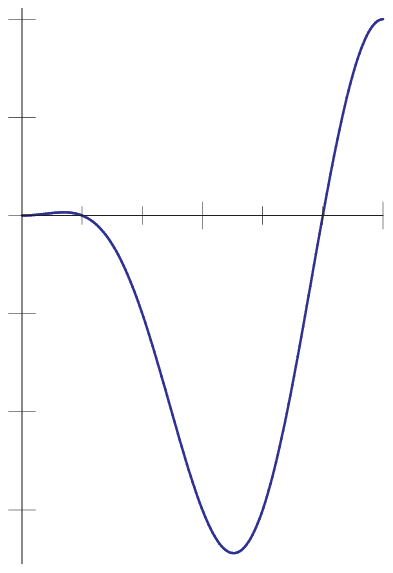}}
    \put(106,81){$\frac{\pi}{2}$}
    \put( 58,81){$\frac{\pi}{4}$}

    \put(-7, 12){$-3$}
    \put(-7, 38.5){$-2$}
    \put(-7, 65){$-1$}
    \put( 3,118){$1$}
    \put( 3,144.5){$2$}
    \put(28,10){$F(\theta)$}
  \end{picture}
% y=-3..2
 \quad
  \begin{picture}(125,153)(-5,0)
    \put(0,0){\includegraphics[height=150pt]{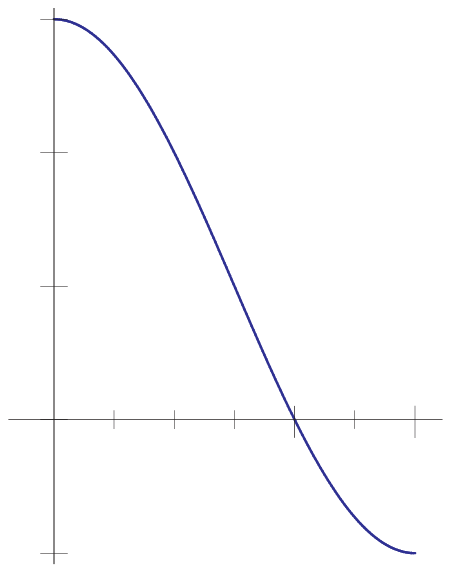}}
    \put(106,27){$\frac{\pi}{2}$}
    \put( 58,27){$\frac{\pi}{4}$}

    \put(-5,  1){$-1$}
    \put( 1, 72){$1$}
    \put( 1,108){$2$}
    \put( 1,143){$3$}
    \put(32,10){$\lambda(\theta)$}
  \end{picture}
% y=-1..3
  \quad
  \begin{picture}(140,153)(-20,0)
    \put(0,0){\includegraphics[height=150pt]{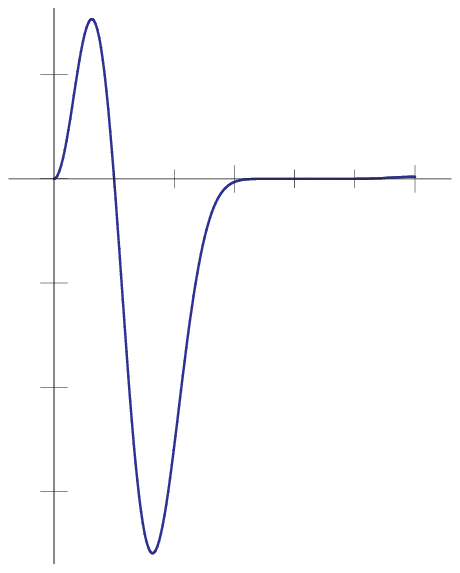}}
    \put(106,92){$\frac{\pi}{2}$}
    \put( 58,92){$\frac{\pi}{4}$}

    \put(-19, 18){$-300$}
    \put(-19, 46){$-200$}
    \put(-19, 74){$-100$}
    \put(-10,130){$100$}
    \put(65,10){$\lambda(\theta)^8 F(\theta)$}
  \end{picture}
% y=-300..100

\caption{The functions $F(\theta)$, $\lambda(\theta)$, and $\lambda(\theta)^8F(\theta)$.}
\label{F:graphs_a=2}
\end{figure}

By the computations in Table~\ref{T:inequality_a=2}, we only need to show that
$K(2^m,4)-K(2^m,1,1)>0$ for $m\geq 14$.
Using~\eqref{Eq:Kostka_integral}, we have
 \begin{multline*}
 \qquad K(2^m,4)-K(2^m,1,1)\ =\ \frac{2}{\pi}\int_0^\pi \lambda_2(\theta)^m
          \bigl( \lambda_4(\theta)\ -\ \lambda_1(\theta)^2\bigr)\; \sin^2{\theta}\;
          d\theta\\
  =\ \frac{2}{\pi}\int_0^\pi  \lambda_2(\theta)^m 
          \bigl( \sin{5\theta}\ \sin{\theta}\ -\ \sin^2{2\theta}\bigr)\; d\theta\\
  =\ \frac{2}{\pi}\int_0^\pi \lambda_2(\theta)^m \frac{1}{2}
         \bigl( 2\cos{4\theta} - \cos{6\theta} -1\bigr)\; d\theta\\
  =\ \frac{1}{\pi}\int_0^\pi \lambda_2(\theta)^m 
         \bigl( 2\cos{4\theta} - \cos{6\theta} -1\bigr)\; d\theta\ .\qquad
 \end{multline*}
The integrand $f(\theta)$ of the last integral is symmetric about $\theta=\pi/2$ in that 
$f(\theta)=f(\pi-\theta)$.
Thus, it  suffices to prove that if $m\geq 14$, then
 \begin{equation}\label{Eq:integral_a=2}
  \int_0^{\frac{\pi}{2}} \lambda_2(\theta)^m (2\cos{4\theta} - \cos{6\theta} -1)
   \, d\theta
   \ >\ 0\,.
\end{equation}
To simplify our notation, set
 \[
   \defcolor{F(\theta)}\ :=\  2\cos 4\theta - \cos\, 6\theta -1
   \qquad\mbox{and}\qquad 
    \defcolor{\lambda(\theta)} := \lambda_2(\theta)\ =\ 1+2\cos{2\theta}\,.
 \]
We display these functions and the integrand in~\eqref{Eq:integral_a=2} for $m=8$
in Figure~\ref{F:graphs_a=2}.
%%%%%%%%%%%%%%%%%%%%%%%%%%%%%%%%%%%%%%%%%%%%%%%%%%%%%%%%%%%%%%%%%%%%%%%%%%%%
%%%%%%%%%%%%%%%%%%%%%%%%%%%%%%%%%%%%%%%%%%%%%%%%%%%%%%%%%%%%%%%%%%%%%%%%%%%%

%

In the interval $[0,\frac{\pi}{2}]$, the zeroes of $F$ occur at $0$, $\frac{\pi}{12}$, and
$\frac{5\pi}{12}$, and $\lambda$ vanishes at $\frac{\pi}{3}$. 
Both functions are positive on $[0,\frac{\pi}{12}]$, and so
 \begin{equation}\label{E:first_reduction}
  \int_{0}^{\frac{\pi}{2}} \lambda^m(\theta) F(\theta)\, d\theta \ \geq \ 
  \int_{0}^{\frac{\pi}{12}} \lambda^m(\theta) F(\theta)\, d\theta\ - \ 
    \int_{\frac{\pi}{12}}^{\frac{\pi}{2}} \big| \lambda^m(\theta) F(\theta)\, \big| \,
    d\theta \,.
 \end{equation}
We show the positivity of \eqref{Eq:integral_a=2} by showing that the right hand side
of~\eqref{E:first_reduction} is positive for $m\geq 14$.
This is equivalent to the following inequality,
\begin{equation}\label{Eq:main_inequality_a=2}
  \int_{0}^{\frac{\pi}{12}} \lambda^m(\theta) F(\theta)\, d\theta\ >\  
  \int_{\frac{\pi}{12}}^{\frac{\pi}{3}} \big| \lambda^m(\theta) F(\theta)\,\big| 
  \, d\theta \ + \ 
  \int_{\frac{\pi}{3}}^{\frac{\pi}{2}} \big| \lambda^m(\theta) F(\theta)\, \big| \, d\theta\,.
\end{equation}

The function $\lambda(\theta)$ is monotone decreasing in the interval $[0,\frac{\pi}{2}]$,
and it vanishes at $\frac{\pi}{3}$, so the maximum of $| \lambda(\theta) | $ on this
interval is $|\lambda(\frac{\pi}{2})| = 1$. 
Also, $|F(\theta)| \leq 4$ for all $\theta \in [0,\frac{\pi}{2}]$.
Thus we estimate the last integral in~\eqref{Eq:main_inequality_a=2},
\[ 
  \int_{\frac{\pi}{3}}^{\frac{\pi}{2}} \big| \lambda^m(\theta) F(\theta) \big|\, d\theta \ \leq \
  \int_{\frac{\pi}{3}}^{\frac{\pi}{2}} 1\cdot 4\, d\theta \ = \ \frac{2\pi}{3}\,.
\]
It is therefore enough to show that
\begin{equation}\label{Eq:RHS_main_inequality}
 \int_{0}^{\frac{\pi}{12}} \lambda^m(\theta) F(\theta)\, d\theta \ > \  
 \int_{\frac{\pi}{12}}^{\frac{\pi}{3}} \big| \lambda^m(\theta) F(\theta)\,\big| 
 \, d\theta\ +\ \frac{2\pi}{3}\,,
\end{equation}
for $m\geq 14$.
We establish~\eqref{Eq:RHS_main_inequality} by induction on $m\geq 14$.
This inequality holds for $m=14$, as the left hand side is
\[
  \int_{0}^{\frac{\pi}{12}} \lambda^{14}(\theta) F(\theta)\, d\theta \ =\ 
  \frac{69}{4}\pi + \frac{26374}{7}\sqrt{3} + \frac{1679543168}{255255}\ \approx\ 13159.9%2898
\]
whereas the right hand side is 
\[
   \int_{\frac{\pi}{12}}^{\frac{\pi}{3}} \big| \lambda^{14}(\theta) F(\theta)\,\big| 
   \, d\theta + \frac{2\pi}{3} = \frac{63052312}{17017} \sqrt{3} -\frac{613}{12}\pi+\frac{1679543168}{255255}
   \ \approx\  12837.1%06957.
\]
Suppose now that the inequality~\eqref{Eq:RHS_main_inequality} holds for some $m\geq 14$.

As $\lambda(\frac{\pi}{12})=1+\sqrt{3}$ and $\lambda$ is decreasing in $[0,\frac{\pi}{2}]$,
we have $\lambda(\theta)\geq 1+\sqrt{3}$ for $\theta\in[0,\frac{\pi}{12}]$.
Thus
 \begin{equation}\label{Eq:first_simplification}
   \int_{0}^{\frac{\pi}{12}}  \lambda^{m+1}(\theta) F(\theta) \, d\theta \ \geq \ 
   \int_{0}^{\frac{\pi}{12}} \left(1{+}\sqrt{3}\right) \cdot \lambda^{m}(\theta) F(\theta) \, d\theta.
 \end{equation}
Similarly, when $\theta \in [\frac{\pi}{12}, \frac{\pi}{2}]$ we have that
$|\lambda(\theta)| \leq 1{+}\sqrt{3}$, as $\lambda(\frac{\pi}{2})=-1$. 
Therefore,
\begin{equation}\label{Eq:second_simplification}
\int_{\frac{\pi}{12}}^{\frac{\pi}{3}} \big| \lambda^{m+1}(\theta) F(\theta)\, \big|\, d\theta \ \leq \ 
\int_{\frac{\pi}{12}}^{\frac{\pi}{3}} \left(1{+}\sqrt{3}\right) \cdot \big| \lambda^{m}(\theta) F(\theta)\, \big|\, d\theta.
\end{equation}

From the induction hypothesis and equations \eqref{Eq:first_simplification} , and \eqref{Eq:second_simplification}, we obtain 
\begin{multline}
\qquad \int_{0}^{\frac{\pi}{12}}  \lambda^{m+1}(\theta) F(\theta) \, d\theta \ \geq \ 
\int_{0}^{\frac{\pi}{12}} \left(1{+}\sqrt{3}\right) \cdot \big| \lambda^{m}(\theta) F(\theta)\, \big|\, d\theta \\
> \ \int_{\frac{\pi}{12}}^{\frac{\pi}{3}} \left(1{+}\sqrt{3}\right) \cdot \big| \lambda^{m}(\theta) F(\theta)\, \big|\, d\theta + (1{+}\sqrt{3})\cdot \frac{2\pi}{3} \\
> \ \int_{\frac{\pi}{12}}^{\frac{\pi}{3}} \big| \lambda^{m+1}(\theta) F(\theta)\, \big|\, d\theta +  \frac{2\pi}{3}. \   
\qquad
\end{multline}
This completes the proof of Lemma~\ref{L:inequality_a=2}.
%\end{proof}

%%%%%%%%%%%%%%%%%%%%%%%%%%%%%%%%%%%%%%%%%%%%%%%
%%%%%%%%%%%%%%%%%%%%%%%%%%%%%%%%%%%%%%%%%%%%%%%%%%%%%%%%%%%%%%%%%%%%%%%%%%%%%%

\providecommand{\bysame}{\leavevmode\hbox to3em{\hrulefill}\thinspace}
\providecommand{\MR}{\relax\ifhmode\unskip\space\fi MR }
\providecommand{\MRhref}[2]{%
  \href{http://www.ams.org/mathscinet-getitem?mr=#1}{#2}
}
\providecommand{\href}[2]{#2}

\end{document}